\documentclass[preprint,12pt]{elsarticle}

\usepackage{graphicx}

\usepackage[ruled,lined,algonl,boxed]{algorithm2e}

\usepackage{amssymb}
\usepackage{amsthm}
\newtheorem{thm}{Theorem}
\newtheorem{cor}[thm]{Corollary}
\newtheorem{conj}[thm]{Conjecture}
\newtheorem{lemma}[thm]{Lemma}
\newdefinition{rmk}{Remark}
\newproof{pf}{Proof}
\newproof{pot}{Proof of Theorem \ref{thm2}}

\journal{\empty}

\begin{document}

\begin{frontmatter}



\title{The bondage number of graphs on topological surfaces and Teschner's conjecture}


\author[label2]{Andrei Gagarin\corref{cor1}}
\ead{andrei.gagarin@acadiau.ca}
\cortext[cor1]{Corresponding author, fax: (1-902)\,585-1074, phone: (1-902)\,585-1419}
\author[label2a]{Vadim Zverovich}
\ead{vadim.zverovich@uwe.ac.uk}
\address[label2]{Department of Mathematics and Statistics, Acadia University, Wolfville, Nova Scotia, B4P 2R6, Canada}
\address[label2a]{Department of Mathematics and Statistics, University of the West of England, Bristol, BS16 1QY, UK}

\begin{abstract}
The bondage number of a graph is the smallest number of its edges whose removal results in a graph having a larger domination number. 
We provide constant upper bounds for the bondage number of graphs on topological surfaces, 
improve upper bounds for the bondage number in terms of the maximum vertex degree and the orientable and non-orientable genera of the graph, 
and show tight lower bounds for the number of vertices of graphs $2$-cell embeddable on topological surfaces of a given genus. 
Also, we provide stronger upper bounds for graphs with no triangles and graphs with the number of vertices larger than a certain threshold in terms of the graph genera. 
This settles Teschner's Conjecture in positive for almost all graphs.
\end{abstract}


\end{frontmatter}


\section{Introduction}
\label{Intro}
We consider simple finite non-empty graphs. 
For a graph $G$, its vertex and edge sets are denoted, respectively, by $V(G)$ and $E(G)$, $|V(G)|=n$ and $|E(G)|=m$.
We also use the following standard notation: $d(v)$ for the degree of a vertex $v$ in $G$,
$\Delta=\Delta(G)$ for the maximum vertex degree of $G$, 
$\delta=\delta(G)$ for the minimum vertex degree of $G$, and 
$N(v)$ for the neighbourhood of a vertex $v$ in $G$. 

A set $D\subseteq V(G)$ is a {\it dominating set} if every vertex not in $D$ is
adjacent to at least one vertex in $D$. The minimum cardinality of
a dominating set of $G$ is the \emph{domination number}
$\gamma(G)$. Clearly, for any spanning subgraph $H$ of $G$, $\gamma(H)\ge\gamma(G)$. The \emph{bondage number} of $G$, denoted by $b(G)$, is the minimum cardinality of a set of edges $B\subseteq E(G)$ such that $\gamma(G-B)>\gamma(G)$, where $V(G-B)=V(G)$ and $E(G-B)=E(G)\backslash B$. In a sense, the bondage number $b(G)$ measures integrity and reliability of the domination number $\gamma(G)$ with respect to the edge removal from $G$, which may correspond, e.g., to link failures in communication networks.

The bondage number was introduced by Bauer et al. \cite{B83} (see also Fink et al. \cite{F90}). Recently, it has been shown by Hu and Xu \cite{HX2011} that the decision problem for the bondage number is NP-hard. Also, they have conjectured that determining an actual set of edges corresponding to the bondage number is not even an NP-problem, which implies it is important to have any reasonable estimations and bounds on the bondage number in terms of other graph parameters and properties.
Two unsolved classical conjectures for the bondage number of arbitrary and planar graphs are as follows.
\begin{conj} [Teschner \cite{T95}] \label{all}
For any graph $G$, $b(G)\le \frac{3}{2}\Delta(G)$.
\end{conj}

Hartnell and Rall \cite{H94} and Teschner \cite{T96} showed that for the cartesian product $G_n = K_n \times K_n$, $n\ge2$, the bound of Conjecture \ref{all} is sharp, i.e. $b(G_n)= \frac{3}{2}\Delta(G_n)$. Teschner \cite{T95} also proved that Conjecture \ref{all} holds when $\gamma(G) \le 3$.

\begin{conj} [Dunbar et al. \cite{D98}] \label{planar-conj}
If $G$ is a planar graph, then $b(G)\le\Delta(G)+1$.
\end{conj}

Trying to prove Conjecture \ref{planar-conj}, Kang and Yuan \cite{K00} have shown the following, with a simpler topological proof later discovered by Carlson and Develin \cite{C06}. 
\begin{thm} [\cite{K00,C06}] \label{planar}
For any connected planar graph $G$, 
$$b(G)\le \min\{8,\ \Delta(G)+2\}.$$
\end{thm}

This solves Conjecture \ref{planar-conj} when $\Delta(G)\ge7$, and Conjecture \ref{all} for planar graphs with $\Delta(G)\ge 4$. 
Also, it is shown in \cite{C06} that $b(G)\le \Delta(G)+3$ for any connected toroidal graph $G$, which solves Conjecture \ref{all} for toroidal graphs with $\Delta(G)\ge 6$.
In \cite{GZ2011}, we generalize this for any topological surface as follows.
\begin{thm} [\cite{GZ2011}] \label{original-thm}
For a connected graph $G$ of orientable genus $h$ and non-orientable genus $k$, 
\begin{equation}\label{original-ub}
b(G)\le \min\{\Delta(G)+h+2,\ \Delta(G)+k+1\}.
\end{equation}
\end{thm}
Also, in \cite{GZ2011}, we indicate that the upper bound (\ref{original-ub}) can be improved for bigger values of the genera $h$ and $k$ by adjusting the proofs and should be helpful in solving Conjecture \ref{all}. Finally, we state the following general conjecture.
\begin{conj} [\cite{GZ2011}]\label{conj-big}
For a connected graph $G$ of orientable genus $h$ and non-orientable genus $k$, $b(G)\le\min\{c_h,\, c^\prime_k,\, \Delta(G)+o(h),\, \Delta(G)+o(k)\}$, where $c_h$ and $c^\prime_k$ are constants depending, respectively, on the orientable and non-orientable genera of $G$.
\end{conj}

Notice that it is sufficient to consider connected graphs because the bondage number of a disconnected graph $G$ is the minimum of the bondage numbers of its components.

In this paper, we provide constant upper bounds for the bondage number of graphs on topological surfaces, which can be used as the first estimation for the constants $c_h$ and $c'_k$ of Conjecture \ref{conj-big}. Also, we improve upper bounds of Theorem \ref{original-thm}, and show tight lower bounds for the number of vertices of graphs $2$-cell embeddable on topological surfaces of a given genus. Also, we provide stronger upper bounds for graphs with no triangles and graphs with the number of vertices larger than a certain threshold in terms of the genera $h$ and $k$. This provides ideas for improvements of our results in case of some restricted classes of graphs, shows that the bondage number is at most eleven and settles Teschner's Conjecture \ref{all} in positive for almost all graphs.




\section{Graphs on the topological surfaces}
\label{Section2}

The planar graphs are precisely the graphs that can be drawn with no crossing edges on the sphere $S_0$.
A topological surface $S$ can be obtained from the sphere $S_0$ by adding a number of handles or crosscaps. If we add $h$, $h\ge 1$, handles to $S_0$, we obtain an orientable surface $S_h$, which is often referred to as the \emph{$h$-holed torus}. The number $h$ is called the \emph{orientable genus} of $S_h$. If we add $k$, $k\ge 1$, crosscaps to the sphere $S_0$, we obtain a non-orientable surface $N_k$. The number $k$ is called the \emph{non-orientable genus} of $N_k$. Any topological surface is homeomorphically equivalent either to $S_h$ ($h\ge 0$), or to $N_k$ ($k\ge 1$). For example, $S_1$, $N_1$, $N_2$ are the \emph{torus}, the \emph{projective plane}, and the \emph{Klein bottle}, respectively.

A graph $G$ is \emph{embeddable} on a topological surface $S$ if it admits a drawing on the surface with no crossing edges. Such a drawing of $G$ on the surface $S$ is called an \emph{embedding} of $G$ on $S$. Notice that there can be many different embeddings of the same graph $G$ on a particular surface $S$. The embeddings can be distinguished and classified by different properties. The set of faces of a particular embedding of $G$ on $S$ is denoted by $F(G)$, $|F(G)|=f$.

An embedding of $G$ on the surface $S$ is a \emph{$2$-cell embedding} if each face of the embedding is homeomorphic to an open disk. In other words, a $2$-cell embedding is an embedding on $S$ that ``fits" the surface. This is expressed in the Euler's formula (\ref{Euler-formula}) of Theorem \ref{Euler} below. For example, a cycle $C_n$ ($n\ge3$) does not have a $2$-cell embedding on the torus, but it has $2$-cell embeddings on the sphere and the projective plane. Similarly, a planar graph may have $2$-cell and non-$2$-cell embeddings on the torus.
An algorithm to transform a planar $2$-cell embedding into a toroidal $2$-cell embedding, whenever possible, can be found in Gagarin et al. \cite{GKN}, pp.\,358--360. Similar algorithms to transform a $2$-cell embedding of genus $h$, $h\ge 1$, (resp., $k$, $k\ge 1$) into a $2$-cell embedding of genus $h+1$ (resp., $k+1$), whenever possible, can be devised for orientable (resp., non-orientable) surfaces by analogy, with more cases to consider. See also how to transform a planar $2$-cell embedding of a graph with a cycle into a projective-planar $2$-cell embedding in Kocay and Kreher \cite{KK2005}, p.\,364.

The following result is usually known as the (generalized) \emph{Euler's formula}. We state it here in a form similar to Thomassen \cite{T92}.
\begin{thm} [Euler's Formula, \cite{T92}] \label{Euler}
Given a connected graph $G$ with $n$ vertices and $m$ edges $2$-cell embedded on a topological surface $S$,
\begin{equation} \label{Euler-formula}
n - m + f = \chi(S),
\end{equation}
where either $\chi(S)=2-2h$ and $S=S_h$, or $\chi(S)=2-k$ and $S=N_k$, and $f$ is the number of faces of the $2$-cell embedding on $S$.
\end{thm}
Equation (\ref{Euler-formula}) is usually referred to as the \emph{Euler's formula} for an orientable surface $S_h$ of genus $h$, $h\ge 0$, or a non-orientable surface $N_k$ of genus $k$, $k\ge 1$, and the invariant $\chi(S)$ is the \emph{Euler characteristic} of an orientable surface $S=S_h$ or a non-orientable surface $S=N_k$, respectively.

The \emph{orientable genus} of a graph $G$ is the smallest integer $h=h(G)$ such that $G$ admits an embedding on an orientable topological surface $S$ of genus $h$.
The \emph{non-orientable genus} of $G$ is the smallest integer $k=k(G)$ such that $G$ can be embedded on a non-orientable topological surface $S$ of genus $k$. In general, $h(G)\not=k(G)$ (e.g., see \cite{KK2005}, pp.\,367-368), and the embeddings on $S_{h(G)}$ and $N_{k(G)}$ must be $2$-cell embeddings.


\begin{lemma} \label{lemma-verts-surfs} 
Given a graph $G$ $2$-cell embedded on an orientable surface $S_h$ of genus $h$,
\begin{equation}\label{b-orient}
|V(G)|=n\ge\frac{3+\sqrt{16h+1}}{2}>2\sqrt{h}+1,\ \ \ h\ge 1,
\end{equation}
$|V(G)|=n\ge1$ for $h=0$, and, on a non-orientable surface $N_k$ of genus $k$,
\begin{equation}\label{b-non-orient}
|V(G)|=n\ge\frac{3+\sqrt{8k+1}}{2}>\sqrt{2k}+1,\ \ \ k\ge 1.
\end{equation}
\end{lemma}

\begin{pf} From Euler's formula (\ref{Euler-formula}),
$$
n - m + f = \chi(S).
$$
Since $f\ge 1$, $m\le\frac{n(n-1)}{2}$, we have
$$
\chi(S)=n-m+f\ge n-\frac{n(n-1)}{2}+1,
$$
which gives
$$
n^2-3n+2(\chi(S)-1)\ge 0.
$$
Solving the corresponding quadratic equation for $n$,
\begin{equation} \label{inequality1}
n=\frac{3\pm\sqrt{17-8\chi(S)}}{2}.
\end{equation}
Since $n$ is a positive integer, plugging in $\chi(S)=2-2h$ and $\chi(S)=2-k$ into (\ref{inequality1}) gives
$$
n\ge\frac{3+\sqrt{16h+1}}{2}>2\sqrt{h}+1\ \ \ \mathrm{for}\ \ h\ge 1,
$$
$n\ge1$ when $h=0$, and 
$$
n\ge\frac{3+\sqrt{8k+1}}{2}>\sqrt{2k}+1\ \ \ \mathrm{for}\ \ k\ge 1,
$$
respectively.
\qed
\end{pf}

The \emph{maximum orientable} (resp., \emph{non-orientable}) \emph{genus} $h_M(G)$ (resp., $k_M(G)$) of a graph $G$ is the largest integer $h$ (resp., $k$) such that $G$ has a $2$-cell embedding on $S_h$ (resp., $N_k$). The maximum genera of graphs are well-studied parameters (for example, see Huang \cite{Huang1} and Ringel \cite{R77}).
Notice that, if $h(G)$ is the orientable genus of $G$, then $G$ has $2$-cell embeddings on the orientable surfaces of genus $h$, $h(G)\le h\le h_M(G)$. Similarly, $G$ has $2$-cell embeddings on the non-orientable surfaces of genus $k$, $k(G)\le k\le k_M(G)$.

The bounds of Lemma \ref{lemma-verts-surfs} are tight. Euler's formula (\ref{Euler-formula}) implies $h_M(G)\le\lfloor\frac{m-n+1}{2} \rfloor$ and $k_M(G)\le m-n+1$, and $4$-edge connected graphs are known to be \emph{upper-embeddable}, i.e. to have $h_M(G)=\lfloor\frac{m-n+1}{2} \rfloor$ (e.g., see Jungerman \cite{J1978}). Notice that complete graphs $K_n$, $n\ge 5$, are $4$-edge connected, and $h(K_n)=\lceil\frac{(n-3)(n-4)}{12}\rceil$ (e.g., see \cite{BM2010}, p.\,281).
Then,\\ 
for $h=3$, Lemma \ref{lemma-verts-surfs} gives $n\ge 5$, and $K_5$ has $h(K_5)=1$, $h_M(K_5)=3$;\\ 
for $h=5$, Lemma \ref{lemma-verts-surfs} gives $n\ge 6$, and $K_6$ has $h(K_6)=1$, $h_M(K_6)=5$;\\ 
for $h=14$, Lemma \ref{lemma-verts-surfs} gives $n\ge 9$, and $K_9$ has $h(K_9)=3$, $h_M(K_9)=14$;\\ 
for $h=18$, Lemma \ref{lemma-verts-surfs} gives $n\ge 10$, and $K_{10}$ has $h(K_{10})=4$, $h_M(K_{10})=18$;\\ 
etc. In general, for $h=4y^2\pm y$, $y\in\mathbb{Z}^+$, the bound of Lemma \ref{lemma-verts-surfs} is tight. Notice that, from the proof of Lemma \ref{lemma-verts-surfs}, the $2$-cell embeddings of $K_{4y^2\pm y}$ have a unique face on the surface of its maximum genus $h_M(K_{4y^2\pm y})$, $y\in\mathbb{Z}^+$. Similar observations can be easily obtained for the non-orientable surfaces: a connected graph $G$ which is not a tree has $k_M(G)=m-n+1$ (see \cite{R77}).

A \emph{triangle-free graph} $G$ is a graph having no cycles of length $3$. The lower bounds of Lemma \ref{lemma-verts-surfs} can be easily improved for graphs with no triangles as follows.
\begin{cor} \label{cor-verts-surfs} 
Given a triangle-free graph $G$ $2$-cell embedded on an orientable surface $S_h$ of genus $h$, the number of vertices of $G$ is
\begin{equation}\label{cor-b-orient}
n\ge 2(1+\sqrt{2h}),\ \ \ h\ge 1,
\end{equation}
$n\ge1$ when $h=0$, and, on a non-orientable surface $N_k$ of genus $k$,
\begin{equation}\label{cor-b-non-orient}
n\ge 2(1+\sqrt{k}),\ \ \ k\ge 1.
\end{equation}
\end{cor}

\begin{pf}
The number of edges of a triangle-free graph $G$ is $m\le \frac{n^2}{4}$ (e.g., see \cite{BM2010}, p.\,45).
The rest of the proof is similar to that of Lemma \ref{lemma-verts-surfs}.
\qed
\end{pf}

Notice that bipartite graphs are triangle-free, and all the results in this paper for the triangle-free graphs hold for the bipartite graphs as well.


\section{Constant upper bounds for general graphs on the topological surfaces}
\label{Section3}

Hartnell and Rall \cite{HR99} proved the following.

\begin{thm} [Hartnell and Rall \cite{HR99}] \label{thm-HR}
The number of edges of a connected graph $G$ with $n$ vertices and the bondage number $b(G)$ is $m\ge\frac{n}{4}(b(G)+1)$, and the bound is sharp.
\end{thm} 

We use Theorem \ref{thm-HR} to establish the following upper bounds on the bondage number of graphs.

\begin{thm}\label{thm-eleven}
Let $G$ be a connected graph of the orientable genus $h=h(G)$, the non-orientable genus $k=k(G)$, and having $n=|V(G)|$ vertices. Then\\

\noindent (i) $h=0$ ($G$ is planar) or $k=1$ ($G$ is projective-planar) implies $b(G)\le 10$;\\

\noindent (ii) $h\ge 1$ and $n>12(2h-2)$, or $k\ge 2$ and $n>12(k-2)$, imply $b(G)\le 11$;\\

\noindent (iii) $h\ge 2$ and $\frac{3+\sqrt{16h+1}}{2}\le n\le12(2h-2)$ imply $b(G)\le 11+\frac{24(h-1)(3-\sqrt{16h+1})}{1-8h}=11+O(\sqrt{h})$, and $k\ge 3$ and $\frac{3+\sqrt{8k+1}}{2}\le n\le12(k-2)$ imply $b(G)\le 11+\frac{12(k-2)(3-\sqrt{8k+1})}{1-4k}=11+O(\sqrt{k})$.
\end{thm} 

\begin{pf} As a corollary of Euler's formula (\ref{Euler-formula}), for $n\ge 3$, in general,
\begin{equation}\label{edge-number}
m\le 3(n-\chi(S))
\end{equation}
(e.g., see \cite{BM2010}, p.\,279). Then, (\ref{edge-number}) and Theorem \ref{thm-HR} give
$$
\frac{n(b(G)+1)}{4}\le m\le 3(n-\chi(S)),
$$
which implies
\begin{equation}\label{upper-eleven}
b(G)\le 11-\frac{12\chi(S)}{n}.
\end{equation}
Plugging in $\chi(S)=2-2h$ and $\chi(S)=2-k$ into (\ref{upper-eleven}) gives
$$
b(G)\le 11+\frac{12(2h-2)}{n}\ \ \ \ \ \mathrm{and}\ \ \ \ \ \ b(G)\le 11+\frac{12(k-2)}{n},
$$
respectively. The statements of Theorem \ref{thm-eleven} follow directly, applying the bounds (\ref{b-orient}) and (\ref{b-non-orient}) of Lemma \ref{lemma-verts-surfs} to obtain the statement of Theorem \ref{thm-eleven}(iii).
\qed
\end{pf}

Clearly, in the case of planar graphs, Theorem \ref{planar} provides a better upper bound, $b(G)\le c_0\le 8$, than Theorem \ref{thm-eleven}(i).
Since there are no restrictions on the number of vertices in the cases of toroidal ($h=1$), projective-planar ($k=1$), and Klein bottle ($k=2$) graphs in Theorem \ref{thm-eleven}, we have the following general constant upper bounds.

\begin{cor}\label{cor-tor-kl-bot} 
For any connected projective-planar graph $G$, $b(G)\le c'_1\le 10$, and any connected toroidal or Klein bottle graph $G$, $b(G)\le 11$, i.e. $c_1\le 11$ and $c'_2\le 11$.
\end{cor}
The formulae of Theorem \ref{thm-eleven}(iii) provide constant upper bounds for the surfaces of higher genera as follows.
\begin{cor}\label{cor-const-ub} 
For a connected graph $G$ of orientable genus $h=h(G)\ge 2$ and non-orientable genus $k=k(G)\ge 3$, we have 

%
\begin{table}[h!] 
\centerline {\footnotesize 
	\begin{tabular}[t]{|| r || r | r | r | r | r | r | r | r | r | r | r | r | r | r ||}
	\hline
	Orientable genus, $h$ & 2 & 3 & 4 & 5 & 6 & 7& 8 & 9 & 10 & 11 & 12 & 13 & 14 & 15\\
	\hline
	$b(G)\le c_h\le$ & 15 & 19 & 22 & 25 & 28 & 30 & 33 & 35 & 37 & 39 & 41 & 43 & 44 & 46\\
	\hline
	\hline
	Non-orientable genus, $k$ & 3 & 4 & 5 & 6 & 7 & 8 & 9 & 10 & 11 & 12 & 13 & 14 & 15 & 16\\
	\hline
	$b(G)\le c'_k\le$ & 13 & 15 & 17 & 19 & 21 & 22 & 24 & 25 & 27 & 28 & 29 & 30 & 32 & 33\\
	\hline
	\end{tabular} }
\caption{Constant upper bounds for the bondage number of graphs on topological surfaces of higher genera ($h\le 15$ and $k\le 16$).\label{tab-average}}
\end{table}
%
\end{cor}
Clearly, since only the direct arguments with Euler's formulae have been used in Theorem \ref{thm-eleven} and Lemma \ref{lemma-verts-surfs}, the bounds for $N_k$, $k=2h$, coincide with the corresponding upper bounds for $S_h$, $h\ge 1$: the Euler characteristics and corresponding Euler's formulae are the same in this case. However, surfaces $S_h$ and $N_{2h}$, $h\ge1$, are not equivalent, i.e. \emph{non-homeomorphic} (e.g., see \cite{T92}, pp.\,129--130), and the classes of graphs of orientable genus $h$ and non-orientable genus $k=2h$, $h\ge 1$, are quite different. Therefore, we conjecture that refinements of the results of Theorem \ref{thm-eleven} are going to provide different constant upper bounds for the bondage number of graphs embeddable on $S_h$ and $N_{2h}$, $h\ge 1$.


\section{Graphs with no triangles}
\label{Section4}

The triangle-free graphs are exactly the graphs of girth at least $4$. 
Fischermann et al. \cite{F03} have shown the following.

\begin{thm} [Fischermann et al. \cite{F03}] \label{thm-FRV}
A connected planar triangle-free graph $G$ has $b(G)\le 6$.
\end{thm} 

We provide a simple proof of Theorem \ref{thm-FRV} and generalize it to all the other topological surfaces as follows.

\begin{thm}\label{thm-no-triangles}
Let $G$ be a connected triangle-free graph of the orientable genus $h=h(G)$, the non-orientable genus $k=k(G)$, and having $n=|V(G)|$ vertices. Then\\

\noindent (i) $h=0$ ($G$ is planar) or $k=1$ ($G$ is projective-planar) implies $b(G)\le 6$;\\

\noindent (ii) $h\ge 1$ and $n>8(2h-2)$, or $k\ge 2$ and $n>8(k-2)$, imply $b(G)\le 7$;\\

\noindent (iii) $h\ge 2$ and $2(1+\sqrt{2h})\le n\le8(2h-2)$ imply $b(G)\le 7+\frac{8(h-1)}{1+\sqrt{2h}}$, and $k\ge 3$ and $2(1+\sqrt{k})\le n\le8(k-2)$ imply $b(G)\le 7+\frac{4(k-2)}{1+\sqrt{k}}$.
\end{thm}

\begin{pf} In case of triangle-free graphs, $4f\le2m$ and $f\le\frac{m}{2}$. Then, similarly to (\ref{edge-number}), as a corollary to Euler's formula (\ref{Euler-formula}), for $n\ge 3$, in general,
\begin{equation}\label{edge-number-2}
m\le 2(n-\chi(S))
\end{equation}
Then, in this case, (\ref{edge-number-2}) and Theorem \ref{thm-HR} give
$$
\frac{n(b(G)+1)}{4}\le m\le 2(n-\chi(S)),
$$
which implies
\begin{equation}\label{upper-seven}
b(G)\le 7-\frac{8\chi(S)}{n}.
\end{equation}
Then, plugging in $\chi(S)=2-2h$ and $\chi(S)=2-k$ into (\ref{upper-seven}) gives
$$
b(G)\le 7+\frac{8(2h-2)}{n}\ \ \ \ \ \mathrm{and}\ \ \ \ \ \ b(G)\le 7+\frac{8(k-2)}{n},
$$
respectively. The statements of Theorem \ref{thm-no-triangles} follow directly, applying the bounds (\ref{cor-b-orient}) and (\ref{cor-b-non-orient}) of Corollary \ref{cor-verts-surfs} to Lemma \ref{lemma-verts-surfs} to obtain the statement of Theorem \ref{thm-no-triangles}(iii).
\qed
\end{pf}

Notice that, in the case of planar graphs, Theorem \ref{thm-no-triangles}(i) provides the same upper bound, $b(G)\le 6$, as the previously known upper bound of Theorem \ref{thm-FRV}. Conclusions similar to Corollaries \ref{cor-tor-kl-bot} and \ref{cor-const-ub} with constant upper bounds for the bondage number of triangle-free graphs on topological surfaces can be drawn from Theorem \ref{thm-no-triangles} as well. In general, Theorem \ref{thm-no-triangles} provides stronger bounds than Theorem \ref{thm-eleven} in the case of triangle-free graphs.


\section{Improved upper bounds in terms of the maximum vertex degree and the genera}
\label{Section5}
In \cite{GZ2011}, we indicate that the upper bound (\ref{original-ub}) can be improved for bigger values of the genera $h$ and $k$ by adjusting the proofs of the corresponding theorems. 
This is  because the assumptions $\delta(G)\ge h+4$, $\delta(G)\ge k+3$, corresponding arguments and reasoning, used respectively in the proofs of Theorems $8$ and $9$ in \cite{GZ2011}, become vacuous for bigger values of the genera $h$ and $k$ in view of the natural upper bounds on the minimum vertex degree:
$\delta(G)\le\lfloor\frac{5+\sqrt{1+48h}}{2}\rfloor$ for $h\ge 1$, 
$\delta(G)\le\lfloor\frac{5+\sqrt{1+24k}}{2}\rfloor$ for $k\ge 2$ (e.g., see Sachs \cite{Sachs}), and $\delta(G)\le 5$ for a planar or projective-planar graph, i.e. when $h=0$ or $k=1$.
The suggested improvements are of the form 
\begin{equation}\label{original-improved}
b(G)\le \min\{\Delta(G)+h-a,\ \Delta(G)+k-b\},
\end{equation}
where $a\ge -1,\ b\ge 0$, $a,b\in\mathbb{Z}$, $h\ge t_a,\ k\ge t_b$, $t_a,t_b\in\mathbb{Z}^+$, and the thresholds $t_a$ and $t_b$ can be found only after going through the corresponding adjusted proofs. The upper bounds obtained by the adjustments are better than the following upper bound which is asymptotically better than the original bound (\ref{original-ub}):
\begin{equation}\label{cor-ub}
b(G)\le \min\{\Delta(G)+\lfloor\frac{3+\sqrt{1+48h}}{2}\rfloor,\ \Delta(G)+\lfloor\frac{3+\sqrt{1+24k}}{2}\rfloor\},
\end{equation}
$h\ge 1$ and $k\ge 1$. Bound (\ref{cor-ub}) can be seen as a simple corollary to Lemma \ref{lem} below and the upper bounds on the minimum vertex degree in terms of the graph genera. 

First, let us show where bounds (\ref{original-improved}) can be at least as good as bound (\ref{cor-ub}). In other words, having fixed $a\ge -1$ and $b\ge 0$, $a,b\in\mathbb{Z}$, for which values of $h$ and $k$ can we have $h-a\le\lfloor\frac{3+\sqrt{1+48h}}{2}\rfloor$ and $k-b\le\lfloor\frac{3+\sqrt{1+24k}}{2}\rfloor$, respectively? Solving corresponding quadratic inequalities gives $h\in[a+\frac{15}{2}-\frac{\sqrt{48a+217}}{2},a+\frac{15}{2}+\frac{\sqrt{48a+217}}{2}]$ and $k\in[b+\frac{9}{2}-\frac{\sqrt{24b+73}}{2},b+\frac{9}{2}+\frac{\sqrt{24b+73}}{2}]$. Then, clearly, thresholds $t_a$ and $t_b$ computed in the adjusted proofs of Theorems $8$ and $9$ in \cite{GZ2011}, resp., must be in the same intervals to guarantee the bounds of the form (\ref{original-improved}) make sense. For example, when $b=5$ in (\ref{original-improved}), the genus $k\in[3,16]$ guarantee that (\ref{original-improved}) is at least as good as (\ref{cor-ub}), and the threshold is $t_b=13$, i.e. $b(G)\le \Delta(G)+k-5$ for $k\ge 13$.

Apparently, the implicit improvements of the form (\ref{original-improved}) suggested in \cite{GZ2011} (or, maybe, even better) are explicitly obtained and presented in Huang \cite{H2011}. In this section, we improve the upper bounds of \cite{GZ2011} and \cite{H2011} as follows.

One of the classical upper bounds on the bondage number can be stated as follows.
\begin{lemma} [Hartnell and Rall \cite{H94}] \label{lem}
For any edge $uv$ in a graph $G$, we have
$b(G)\le d(u)+d(v)-1-d_{uv}$, where $d_{uv}=|N(u)\cap N(v)|$. 
In particular, this implies that
$b(G) \le \Delta(G)+\delta(G)-1$ \emph{(see also \cite{B83,F90})}.
\end{lemma}

Having a graph $G$ embedded on a surface $S$, each edge $e_i=uv \in E(G)$, $i=1,\ldots,m$, is assigned two weights, $w_i=\frac{1}{d(u)}+\frac{1}{d(v)}$ and $f_i=\frac{1}{m'}+\frac{1}{m''}$, where $m'$ is the number of edges on the boundary of a face on one side of $e_i$, and $m''$ is the number of edges on the boundary of the face on the other side of $e_i$. Notice that, in an embedding on a surface, an edge $e_i$ may be not separating two distinct faces, but instead can appear twice on the boundary of the same face, and, in this case, 
$f_i=\frac{2}{m'}=\frac{2}{m''}$. 

We have 
$$\sum_{i=1}^{m}w_i=n,\ \ \ \ \ \sum_{i=1}^{m}f_i=f.$$ 
and, by Euler's formula (\ref{Euler-formula}), 
$$\sum_{i=1}^{m}(w_i+f_i-1)=n+f-m=\chi(S),$$
or, in other words,
\begin{equation}\label{zero-sum}
\sum_{i=1}^{m}\left(w_i+f_i-1-\frac{\chi(S)}{m}\right)=0.
\end{equation}
Now, each edge $e_i=uv \in E(G)$, $i=1,\ldots,m$, can be associated with the weight $Q(e_i)=w_i+f_i-1-\frac{\chi(S)}{m}$ called, depending on $S$, the \emph{oriented} or \emph{non-oriented} \emph{curvature} of the edge, respectively.

\begin{thm} \label{thm-orient}
Let $G$ be a connected graph $2$-cell embeddable on an orientable surface of genus $h\ge 0$. Then
\begin{equation} \label{eqn-orient}
b(G)\>\le\>\left\{
\begin{array}{ll}
\Delta(G)+\lceil h^{0.7}\rceil+2,&\mathrm{for}\ \ h\le 5,\\
\Delta(G)+\lceil h^{0.7}\rceil+3,&\mathrm{for}\ \ h\ge 6.
\end{array}\right.
\end{equation}
\end{thm}

\begin{pf} 
Suppose $G$ is $2$-cell embedded on the $h$-holed torus $S_h$. Denote by 
$$
\tau\>=\>\left\{
\begin{array}{ll}
\lceil h^{0.7}\rceil -1,&\mathrm{for}\ \ h\le 5,\\
\lceil h^{0.7}\rceil,&\mathrm{for}\ \ h\ge 6.
\end{array}\right.
$$
Then we have to prove
$$
b(G)\le \Delta(G)+\tau+3.
$$
If $\delta(G)\le\tau+4$, then, by Lemma \ref{lem}, 
$$
b(G)\le\Delta(G)+\delta(G)-1\le\Delta(G)+\tau+3,
$$
as required, and inequality (\ref{eqn-orient}) holds. 

Therefore, we can assume $\delta(G)\ge \tau+5$. 
Suppose the opposite, $b(G)\ge \Delta(G)+\tau+4$. Then, by Lemma \ref{lem}, for any edge $e_i=uv$, $i=1,\dots,m$, we have
$$d(u)+d(v)-1-d_{uv}\ge b(G)\ge \Delta(G)+\tau+4,$$
which implies
\begin{equation}\label{star}
d(u)+d(v)\ge\Delta(G)+\tau+5+d_{uv},
\end{equation}
and $\Delta(G)\ge d(u)$, $\Delta(G)\ge d(v)$.
Without loss of generality, assume $d(u)\le d(v)$. Then we have the following three cases.\\ 
 
\emph{Case 1:} $d(u)=\tau+5$. 
By (\ref{star}), $d(v)\ge \Delta(G)+d_{uv}$, which implies $d(v)=\Delta(G)$ and $d_{uv}=0$. Therefore, $m'\ge 4$ and $m''\ge 4$ in this case, and
$Q(e_i)=w_i+f_i-1+\frac{2h-2}{m}$. When $h=0$, $\tau=-1$, and 
$$Q(e_i)\le \frac{2}{4}-\frac{1}{2}-\frac{2}{m}<0.$$
When $h\ge 1$, $d_{uv}=0$ implies $n\ge d(u)+d(v)\ge 2\tau+10$, and $m\ge\frac{n\delta(G)}{2}\ge\frac{(2\tau+10)(\tau+5)}{2}=(\tau+5)^2$, and
$$Q(e_i)\le \frac{2}{\tau+5}-\frac{1}{2}+\frac{2h-2}{(\tau+5)^2}<0.$$

\emph{Case 2:} $d(u)=\tau+6$.
Then $d(v)\ge \tau+6$.
By (\ref{star}), $d(v)\ge \Delta(G)-1+d_{uv}$. If $d_{uv}\ge 2$, then $d(v)\ge \Delta(G)+1$, a contradiction. Therefore, $d_{uv}\le 1$ and, without loss of generality, $m'\ge 3$ and $m''\ge 4$ in this case. We have
$Q(e_i)\le \frac{2}{\tau+6}+\frac{1}{3}+\frac{1}{4}-1+\frac{2h-2}{m}=\frac{2}{\tau+6}-\frac{5}{12}+\frac{2h-2}{m}$.
When $h=0$, $\tau=-1$, and 
$$Q(e_i)\le \frac{2}{5}-\frac{5}{12}-\frac{2}{m}=-\frac{1}{60}-\frac{2}{m}<0.$$
When $h\ge 1$, $d_{uv}\le 1$ implies $n\ge d(u)+d(v)-1\ge 2\tau+11$, and $m\ge\frac{d(u)+d(v)+(n-2)\delta(G)}{2}\ge\frac{2\tau+12+(2\tau+9)(\tau+5)}{2}=\frac{2\tau^2+21\tau+57}{2}$, and
$$Q(e_i)\le \frac{2}{\tau+6}-\frac{5}{12}+\frac{2(2h-2)}{2\tau^2+21\tau+57}<0.$$

\emph{Case 3:} $d(u)\ge\tau+7$.
Then $d(v)\ge \tau+7$, and $m'\ge 3$, $m''\ge 3$. 
We have
$Q(e_i)\le \frac{2}{\tau+7}-\frac{1}{3}+\frac{2h-2}{m}$.
When $h=0$, $\tau=-1$, and 
$$Q(e_i)\le \frac{2}{6}-\frac{1}{3}-\frac{2}{m}<0.$$
When $h\ge 1$, inequality (\ref{star}) implies $n\ge d(u)+d(v)-d_{uv}\ge \Delta(G)+\tau+5\ge 2\tau+12$, then
$m\ge\frac{d(u)+d(v)+(n-2)\delta(G)}{2}\ge\frac{2\tau+14+(2\tau+10)(\tau+5)}{2}=\tau+7+(\tau+5)^2=\tau^2+11\tau+32$, and
$$Q(e_i)\le \frac{2}{\tau+7}-\frac{1}{3}+\frac{2h-2}{\tau^2+11\tau+32}<0\ \ \ \ \mathrm{for}\ \ h\in\mathbb{Z}^+.$$

Thus, $Q(e_i)<0$, $i=1,\ldots,m$, for any $h\ge 0$, $h\in\mathbb{Z}$, and 
$\sum_{i=1}^{m}Q(e_i)<0$, a contradiction with (\ref{zero-sum}). Therefore $b(G)\le\Delta(G)+\tau+3$, as required.
\qed
\end{pf}

If the genus $h$ is fixed, by excluding a finite number of graphs, we can improve the bound of Theorem \ref{thm-orient} as follows.

\begin{cor}\label{cor-orient}
For a connected graph $G$ $2$-cell embeddable on an orientable surface of genus $h\ge 1$, we have:\\
(a) $b(G)\le\Delta(G)+\lceil\ln^2 h\rceil+3$ if $n\ge h$;\\
(b) $b(G)\le\Delta(G)+\lceil\ln h\rceil+3$ if $n\ge h^{1.9}$;\\
(c) $b(G)\le\Delta(G)+4$ if $n\ge h^{2.5}$.
\end{cor}

\begin{pf} 
Denote by:\\
(a) $\tau=\lceil\ln^2 h\rceil$;\\
(b) $\tau=\lceil\ln h\rceil$;\\
(c) $\tau=1$.\\
Then, as in Theorem \ref{thm-orient}, we have to prove
$$
b(G)\le \Delta(G)+\tau+3.
$$
If $\delta(G)\le\tau+4$, then, by Lemma \ref{lem}, 
$$
b(G)\le\Delta(G)+\delta(G)-1\le\Delta(G)+\tau+3,
$$
as required.

Therefore, as in Theorem \ref{thm-orient}, we can assume $\delta(G)\ge \tau+5$,
suppose the opposite, $b(G)\ge \Delta(G)+\tau+4$, and, by Lemma \ref{lem}, for any edge $e_i=uv$, $i=1,\dots,m$, the inequality (\ref{star}) holds here as well.
Again, without loss of generality, $d(u)\le d(v)$. Then we have the following three cases.\\
 

\emph{Case 1:} $d(u)=\tau+5$. 
By (\ref{star}), $d(v)\ge \Delta(G)+d_{uv}$, $d(v)=\Delta(G)$, and $d_{uv}=0$. Therefore, $m'\ge 4$ and $m''\ge 4$. Then
$Q(e_i)=w_i+f_i-1+\frac{2h-2}{m}\le \frac{2}{\tau+5}-\frac{1}{2}+\frac{2h-2}{m}$, and
$m\ge\frac{n\delta(G)}{2}\ge\frac{h^x(\tau+5)}{2}$ and, similarly to the proof of Theorem \ref{thm-orient},
$m\ge\frac{n\delta(G)}{2}\ge\frac{(2\tau+10)(\tau+5)}{2}=(\tau+5)^2$. This implies, respectively, 
$$Q(e_i)=Q'(e_i)\le \frac{2}{\tau+5}-\frac{1}{2}+\frac{2(2h-2)}{h^x(\tau+5)}$$
and
$$Q(e_i)=Q''(e_i)\le \frac{2}{\tau+5}-\frac{1}{2}+\frac{2h-2}{(\tau+5)^2}.$$
We have:\\
(a) $x=1$, then $Q'(e_i)<0$ if $h\ge 12$, and $Q''(e_i)<0$ if $h\le 11$;\\
(b) $x=1.9$, then $Q'(e_i)<0$ if $h\not= 2$, and $Q''(e_i)<0$ if $h=2$;\\
(c) $x=2.5$, then $Q'(e_i)<0$, $h\ge 1$.\\


\emph{Case 2:} $d(u)=\tau+6$.
Then, similarly to the proof of Theorem \ref{thm-orient}, $d(v)\ge \tau+6$, and, by (\ref{star}), $d(v)\ge \Delta(G)-1+d_{uv}$, $d_{uv}\le 1$, and, without loss of generality, $m'\ge 3$ and $m''\ge 4$. 
Then
$Q(e_i)\le \frac{2}{\tau+6}-\frac{5}{12}+\frac{2h-2}{m}$, and $m\ge\frac{d(u)+d(v)+(n-2)\delta(G)}{2}\ge\frac{2\tau+12+(n-2)(\tau+5)}{2}=\frac{2+n(\tau+5)}{2}\ge\frac{2+h^x(\tau+5)}{2}$ and,
similarly to the proof of Theorem \ref{thm-orient}, 
$m\ge\frac{d(u)+d(v)+(n-2)\delta(G)}{2}\ge\frac{2\tau+12+(2\tau+9)(\tau+5)}{2}=\frac{2\tau^2+21\tau+57}{2}.$ This implies, respectively, 
$$Q(e_i)=Q'(e_i)\le \frac{2}{\tau+6}-\frac{5}{12}+\frac{2(2h-2)}{2+h^x(\tau+5)}$$
and
$$Q(e_i)=Q''(e_i)\le \frac{2}{\tau+6}-\frac{5}{12}+\frac{2(2h-2)}{2\tau^2+21\tau+57}.$$
We have:\\
(a) $x=1$, then $Q'(e_i)<0$ if $h\ge 17$, and $Q''(e_i)<0$ if $h\le 16$;\\
(b) $x=1.9$, then $Q'(e_i)<0$ if $h\not= 2$, and $Q''(e_i)<0$ if $h=2$;\\
(c) $x=2.5$, then $Q'(e_i)<0$, $h\ge 1$.\\


\emph{Case 3:} $d(u)\ge\tau+7$.
Then we have $d(v)\ge \tau+7$, $m'\ge 3$, $m''\ge 3$,
$Q(e_i)\le \frac{2}{\tau+7}-\frac{1}{3}+\frac{2h-2}{m}$,
and $m\ge\frac{d(u)+d(v)+(n-2)\delta(G)}{2}\ge\frac{2\tau+14+(n-2)(\tau+5)}{2}=\frac{4+n(\tau+5)}{2}\ge\frac{4+h^x(\tau+5)}{2}$ and,
similarly to the proof of Theorem \ref{thm-orient}, 
$m\ge\frac{d(u)+d(v)+(n-2)\delta(G)}{2}\ge\frac{4+n(\tau+5)}{2}\ge\frac{4+(2\tau+12)(\tau+5)}{2}=2+(\tau+6)(\tau+5)=\tau^2+11\tau+32$. This implies, respectively, 
$$Q(e_i)=Q'(e_i)\le \frac{2}{\tau+7}-\frac{1}{3}+\frac{2(2h-2)}{4+h^x(\tau+5)}$$
and
$$Q(e_i)=Q''(e_i)\le \frac{2}{\tau+7}-\frac{1}{3}+\frac{2h-2}{\tau^2+11\tau+32}.$$
We have:\\
(a) $x=1$, then $Q'(e_i)<0$ if $h\ge 28$, and $Q''(e_i)<0$ if $h\le 27$;\\
(b) $x=1.9$, then $Q'(e_i)<0$ if $h\ge 5$, and $Q''(e_i)<0$ if $h\le 4$;\\
(c) $x=2.5$, then $Q'(e_i)<0$ if $h\not=2$, and $Q''(e_i)<0$ if $h=2$.\\

Thus, $Q(e_i)<0$, $i=1,\ldots,m$, for any $h\ge 1$, $h\in\mathbb{Z}$, and 
$\sum_{i=1}^{m}Q(e_i)<0$, a contradiction with (\ref{zero-sum}). Therefore $b(G)\le\Delta(G)+\tau+3$, as required.
\qed
\end{pf}


\begin{thm} \label{thm-non-orient}
Let $G$ be a connected graph $2$-cell embeddable on a non-orientable surface of genus $k\ge 1$. Then
\begin{equation} \label{eqn-non-orient}
b(G)\>\le\>\left\{
\begin{array}{ll}
\Delta(G)+\lceil k^{0.6}\rceil+1,&\mathrm{for}\ \ k\le 5,\\
\Delta(G)+\lceil k^{0.6}\rceil+2,&\mathrm{for}\ \ k\ge 6.
\end{array}\right.
\end{equation}
\end{thm}

\begin{pf} 
The proof is similar to that of Theorem \ref{thm-orient} above and goes as follows.
Let $G$ be $2$-cell embedded on $N_k$. Denote by 
$$
\tau\>=\>\left\{
\begin{array}{ll}
\lceil k^{0.6}\rceil -1,&\mathrm{for}\ \ k\le 5,\\
\lceil k^{0.6}\rceil,&\mathrm{for}\ \ k\ge 6.
\end{array}\right.
$$
Then we have to prove that
$$
b(G)\le \Delta(G)+\tau+2.
$$
If $\delta(G)\le\tau+3$, then, by Lemma \ref{lem}, 
$$
b(G)\le\Delta(G)+\delta(G)-1\le\Delta(G)+\tau+2,
$$
as required, and inequality (\ref{eqn-non-orient}) holds. 

Therefore, assume $\delta(G)\ge \tau+4$. 
Suppose that $b(G)\ge \Delta(G)+\tau+3$. Then, by Lemma \ref{lem}, for any edge $e_i=uv$, $i=1,\dots,m$,
$$d(u)+d(v)-1-d_{uv}\ge b(G)\ge \Delta(G)+\tau+3,$$
i.e.
\begin{equation}\label{star2}
d(u)+d(v)\ge\Delta(G)+\tau+4+d_{uv}.
\end{equation}
Without loss of generality, $d(u)\le d(v)$. We have the following three cases.\\
 
\emph{Case 1:} $d(u)=\tau+4$. 
By (\ref{star2}), $d(v)\ge \Delta(G)+d_{uv}$, which implies $d(v)=\Delta(G)$, $d_{uv}=0$, and $m'\ge 4$, $m''\ge 4$. In this case, we have
$Q(e_i)=w_i+f_i-1+\frac{k-2}{m}\le\frac{2}{\tau+4}-\frac{1}{2}+\frac{k-2}{m}$. When $k=1$, 
$$Q(e_i)\le \frac{2}{4}-\frac{1}{2}-\frac{1}{m}<0.$$
When $k\ge 2$, $d_{uv}=0$ implies $n\ge d(u)+d(v)\ge 2\tau+8$, and $m\ge\frac{n\delta(G)}{2}\ge\frac{(2\tau+8)(\tau+4)}{2}=(\tau+4)^2$, and
$$Q(e_i)\le \frac{2}{\tau+4}-\frac{1}{2}+\frac{k-2}{(\tau+4)^2}<0.$$

\emph{Case 2:} $d(u)=\tau+5$.
Then $d(v)\ge \tau+5$.
By (\ref{star2}), $d(v)\ge \Delta(G)-1+d_{uv}$, which implies $d_{uv}\le 1$, and, without loss of generality, $m'\ge 3$ and $m''\ge 4$ in this case. Then we have
$Q(e_i)\le \frac{2}{\tau+5}-\frac{5}{12}+\frac{k-2}{m}$.
When $k=1$,
$$Q(e_i)\le \frac{2}{5}-\frac{5}{12}-\frac{1}{m}=-\frac{1}{60}-\frac{1}{m}<0.$$
When $k\ge 2$, $n\ge d(u)+d(v)-1\ge 2\tau+9$, and $m\ge\frac{d(u)+d(v)+(n-2)\delta(G)}{2}\ge\frac{2\tau+10+(2\tau+7)(\tau+4)}{2}=\frac{2\tau^2+17\tau+38}{2}=\tau^2+8.5\tau+19$, and
$$Q(e_i)\le \frac{2}{\tau+5}-\frac{5}{12}+\frac{k-2}{\tau^2+8.5\tau+19}<0.$$

\emph{Case 3:} $d(u)\ge\tau+6$.
Then $d(v)\ge \tau+6$, and $m'\ge 3$, $m''\ge 3$. 
We have
$Q(e_i)\le \frac{2}{\tau+6}-\frac{1}{3}+\frac{k-2}{m}$.
When $k=1$,
$$Q(e_i)\le -\frac{1}{m}<0.$$
When $k\ge 2$, inequality (\ref{star2}) implies $n\ge d(u)+d(v)-d_{uv}\ge \Delta(G)+\tau+4\ge 2\tau+10$. Then
$m\ge\frac{d(u)+d(v)+(n-2)\delta(G)}{2}\ge\frac{2\tau+12+(2\tau+8)(\tau+4)}{2}=\tau+6+(\tau+4)^2=\tau^2+9\tau+22$, and
$$Q(e_i)\le \frac{2}{\tau+6}-\frac{1}{3}+\frac{k-2}{\tau^2+9\tau+22}<0\ \ \ \ \mathrm{for}\ \ k\in\mathbb{Z}.$$

Thus, $Q(e_i)<0$, $i=1,\ldots,m$, for any $k\ge 1$, $k\in\mathbb{Z}$, and 
$\sum_{i=1}^{m}Q(e_i)<0$, a contradiction with (\ref{zero-sum}). Therefore $b(G)\le\Delta(G)+\tau+2$, as required.
\qed
\end{pf}


The bounds of Theorem \ref{thm-non-orient} can be improved if the number of vertices $n$ is restricted from below by some function of the non-orientable genus $k$. If $k$ is fixed, then only a finite number of graphs is forbidden, i.e. the following bounds are true for almost all graphs.

\begin{cor}\label{cor-non-orient}
Let $G$ be a connected graph $2$-cell embeddable on a non-orientable surface of genus $k\ge 1$. Then:\\
(a) $b(G)\le\Delta(G)+\lceil\ln^2 k\rceil+2$ if $n\ge k/6$;\\
(b) $b(G)\le\Delta(G)+\lceil\ln k\rceil+2$ if $n\ge k^{1.6}$;\\
(c) $b(G)\le\Delta(G)+3$ if $n> k^{2}$.
\end{cor}

\begin{pf} For $k=1$, the result follows from Theorem \ref{thm-non-orient}. Therefore, consider $k\ge 2$. Denote by:\\
(a) $\tau=\lceil\ln^2 k\rceil$;\\
(b) $\tau=\lceil\ln k\rceil$;\\
(c) $\tau=1$.\\
In each case, similarly to the proof of Theorem \ref{thm-non-orient}, we have to prove
$$
b(G)\le \Delta(G)+\tau+2.
$$

Again, similarly to the proof of Theorem \ref{thm-non-orient}, we can assume $\delta(G)\ge \tau+4$,
and suppose the opposite, $b(G)\ge \Delta(G)+\tau+3$. Then, by Lemma \ref{lem}, for any edge $e_i=uv$, $i=1,\dots,m$, the inequality (\ref{star2}) holds in this new context as well.
Again, without loss of generality, $d(u)\le d(v)$.\\
 

\emph{Case 1:} $d(v)\ge d(u)=\tau+4$. 
We have
$Q(e_i)\le \frac{2}{\tau+4}-\frac{1}{2}+\frac{k-2}{m}$.\\
(a) $m\ge\frac{n\delta(G)}{2}\ge\frac{k(\tau+4)}{12}$, then $Q(e_i)\le \frac{2}{\tau+4}-\frac{1}{2}+\frac{12(k-2)}{k(\tau+4)}<0$ if $k\ge 122$, and,
similarly to the proof of Theorem \ref{thm-non-orient}, 
$Q(e_i)\le \frac{2}{\tau+4}-\frac{1}{2}+\frac{k-2}{(\tau+4)^2}<0$ if $k\le 121$.\\
(b) $m\ge\frac{n\delta(G)}{2}\ge\frac{k^{1.6}(\tau+4)}{2}$, then $Q(e_i)\le \frac{2}{\tau+4}-\frac{1}{2}+\frac{2(k-2)}{k^{1.6}(\tau+4)}<0$ for $k\ge 2$.\\
(c) $m\ge\frac{n\delta(G)}{2}>\frac{k^2(\tau+4)}{2}$, then $Q(e_i)<\frac{2}{\tau+4}-\frac{1}{2}+\frac{2(k-2)}{k^2(\tau+4)}<0$ for $k\ge 2$.\\


\emph{Case 2:} $d(v)\ge d(u)=\tau+5$,
$Q(e_i)\le \frac{2}{\tau+5}-\frac{5}{12}+\frac{k-2}{m}$.\\
(a) $m\ge\frac{k(\tau+4)}{12}$, then $Q(e_i)\le \frac{2}{\tau+5}-\frac{5}{12}+\frac{12(k-2)}{k(\tau+4)}<0$ if $k\ge 219$, and,
similarly to the proof of Theorem \ref{thm-non-orient}, 
$Q(e_i)\le \frac{2}{\tau+5}-\frac{5}{12}+\frac{k-2}{\tau^2+8.5\tau+19}<0$ if $k\le 218$.\\
(b) $m\ge\frac{d(u)+d(v)+(n-2)\delta(G)}{2}\ge\frac{2\tau+10+(k^{1.6}-2)(\tau+4)}{2}=\frac{k^{1.6}(\tau+4)+2}{2}$, then $Q(e_i)\le \frac{2}{\tau+5}-\frac{5}{12}+\frac{2(k-2)}{k^{1.6}(\tau+4)+2}<0$ for $k\ge 2$.\\
(c) $m\ge\frac{d(u)+d(v)+(n-2)\delta(G)}{2}>\frac{2\tau+10+(k^{2}-2)(\tau+4)}{2}=\frac{k^{2}(\tau+4)+2}{2}$, then $Q(e_i)<\frac{2}{\tau+5}-\frac{5}{12}+\frac{2(k-2)}{k^{2}(\tau+4)+2}<0$ for $k\ge 2$.\\


\emph{Case 3:} $d(v)\ge d(u)=\tau+6$,
$Q(e_i)\le \frac{2}{\tau+6}-\frac{1}{3}+\frac{k-2}{m}$.\\
(a) $m\ge\frac{k(\tau+4)}{12}$, then $Q(e_i)\le \frac{2}{\tau+6}-\frac{1}{3}+\frac{12(k-2)}{k(\tau+4)}<0$ if $k\ge 439$, and,
similarly to the proof of Theorem \ref{thm-non-orient}, 
$Q(e_i)\le \frac{2}{\tau+6}-\frac{1}{3}+\frac{k-2}{\tau^2+9\tau+22}<0$ if $k\le 438$.\\
(b) $m\ge\frac{d(u)+d(v)+(n-2)\delta(G)}{2}\ge\frac{2\tau+12+(k^{1.6}-2)(\tau+4)}{2}=\frac{k^{1.6}(\tau+4)+4}{2}$, then $Q(e_i)\le \frac{2}{\tau+6}-\frac{1}{3}+\frac{2(k-2)}{k^{1.6}(\tau+4)+4}<0$ for $k\ge 2$, $k\in \mathbb{Z}$.\\
(c) $m\ge\frac{d(u)+d(v)+(n-2)\delta(G)}{2}>\frac{2\tau+12+(k^{2}-2)(\tau+4)}{2}=\frac{k^{2}(\tau+4)+4}{2}$, then $Q(e_i)<\frac{2}{\tau+6}-\frac{1}{3}+\frac{2(k-2)}{k^{2}(\tau+4)+4}<0$ for $k\ge 2$.\\

Thus, $Q(e_i)<0$, $i=1,\ldots,m$, for any $k\ge 2$, $k\in\mathbb{Z}$, and 
$\sum_{i=1}^{m}Q(e_i)<0$, a contradiction with (\ref{zero-sum}). Therefore $b(G)\le\Delta(G)+\tau+2$, as required.
\qed
\end{pf}



\section{Teschner's conjecture and final remarks}
\label{Section6}
Combining Corollary \ref{cor-tor-kl-bot} with the results of \cite{GZ2011} and Theorem $6$ in \cite{C06}, we obtain a statements similar to Theorem \ref{planar}.
\begin{thm}\label{tor-kl-bot}
For any connected projective-planar graph $G$,
$$b(G)\le \min\{10,\ \Delta(G)+2\},$$
and any connected toroidal or Klein bottle graph $G$, 
$$b(G)\le \min\{11,\ \Delta(G)+3\}.$$
\end{thm}

Theorem \ref{tor-kl-bot} and Theorem \ref{planar} settle Teschner's Conjecture \ref{all} in positive for all the planar and projective-planar graphs $G$ with $\Delta(G)\ge 4$, and all the toroidal and Klein bottle graphs $G$ with $\Delta(G)\ge 6$, respectively. Theorem \ref{thm-eleven}(ii), Conjecture \ref{cor-orient}(c), and Conjecture \ref{cor-non-orient}(c) settle Teschner's Conjecture \ref{all} in positive for almost all graphs $G$ with $\Delta(G)\ge 6$ or $\Delta(G)\ge 8$, depending on the number of vertices $n$ in terms of the genera $h$ and $k$ of $G$.

Since, for each particular surface $S$ of orientable genus $h$ or non-orientable genus $k$, the number of graphs embeddable in $S$ is infinite, and the number of graphs not satisfying Theorem \ref{thm-eleven}(ii), Conjecture \ref{cor-orient}(c), Conjecture \ref{cor-non-orient}(c), Theorem \ref{tor-kl-bot} or Theorem \ref{planar}, and having the maximum vertex degree $\Delta(G)$ bounded by $8$, $6$, or $4$, respectively, is finite, we can conclude that Teschner's Conjecture \ref{all} holds for almost all graphs in general. Also, this provides particular cases where Teschner's Conjecture \ref{all} must be verified to be settled for all the graphs in general. We hope the hierarchies of upper bounds by the graph genera are going to be useful to solve Conjecture \ref{all} for all the graphs.

In view of the results of \cite{H2011} and upper bounds (\ref{cor-ub}), the bounds of Theorems \ref{thm-orient} and \ref{thm-non-orient} are asymptotically not tight for larger values of the genera $h=h(G)$ and $k=k(G)$. However, Theorems \ref{thm-orient} and \ref{thm-non-orient} provide better results than the upper bounds of \cite{H2011} for values of the genera $h$ and $k$ raging up to several hundreds. Therefore, we would suggest searching for general upper bounds of the form $\Delta(G)+O(\sqrt{h})$ and $\Delta(G)+O(\sqrt{k})$ by combining or improving the results of Theorems \ref{thm-orient} and \ref{thm-non-orient} with those of \cite{H2011}.

Computations in the proofs of Theorems \ref{thm-orient} and \ref{thm-non-orient} and Corollaries \ref{cor-orient} and \ref{cor-non-orient} are done by using software Maple 9.5 and MS Excel. Also, in general, as presented, these four proofs rely on asumptions of appropriate asymptotic behaviour of functions. It would be interesting to obtain similar proofs without using computational tools or the asumptions. However, since asymptotically the upper bounds (\ref{cor-ub}) and corresponding upper bounds of \cite{H2011} are stronger, one needs to do the computations only for a finite number of cases, and there is no need to justify the asymptotic behaviour of functions in general.


For the constant upper bounds for topological surfaces of Theorems \ref{thm-eleven} and \ref{planar} (see also Corollaries \ref{cor-tor-kl-bot} and \ref{cor-const-ub}), we would suggest to refine them to obtain tight constant upper bounds.
For example, Fischermann et al. \cite{F03} asked whether there exist planar graphs of bondage numbers $6$, $7$, or $8$. A class of planar graphs with the bondage number equal to $6$ is shown in \cite{C06}, and $6\le c_0\le 8$ in the case of planar graphs. The next surfaces to consider should be the torus $S_1$, projective plane $N_1$, and Klein bottle $N_2$. Therefore, we would suggest to try to improve, if possible, the results of Corollary \ref{cor-tor-kl-bot} (see also Theorem \ref{tor-kl-bot}).





\end{document}